\documentclass[12pt]{article}
\usepackage[T2A]{fontenc}
\usepackage[utf8]{inputenc}
\usepackage{amssymb,amsmath,amsfonts,amsthm,amscd,latexsym,indentfirst,verbatim}
\def\Span{\mathrm{Span}}
\def\Lie{\mathrm{Lie}}
\def\As{\mathrm{As}}

\begin{document}

\hfill{17B80 (MSC2010)}

\begin{center}
{\Large
Universal enveloping algebra \\of a pair of compatible Lie brackets}

\smallskip

Vsevolod Gubarev
\end{center}

\begin{abstract}
Applying the Poincar\'{e}---Birkhoff---Witt property 
and the Gr\"{o}bner---Shirshov bases technique, 
we find the linear basis of the associative universal enveloping algebra 
in the sense of V. Ginzburg and M. Kapranov of a pair of compatible Lie brackets.
We state that the growth rate of this universal enveloping over $n$-dimen\-sional compatible Lie algebra equals~$n+1$.

\end{abstract}

\section{Introduction}

Hamiltonian pairs (or bihamiltonian structures) play
an important role~\cite{GelDor,GelZakh,Magri}
in the theory of integrable systems from mathematical physics.
Such structures correspond to pairs of compatible Poisson brackets 
defined on the same manifold. 
Two Poisson brackets $\{\cdot,\cdot\}_1$ and $\{\cdot,\cdot\}_2$ 
are said to be compatible if
$\alpha \{\cdot,\cdot\}_1 + \beta \{\cdot,\cdot\}_2$ 
is a~Poisson bracket for all $\alpha,\beta\in \Bbbk$,
where $\Bbbk$~denotes the ground field.
In terms of operads, algebras with compatible Poisson brackets form 
a so called bi-Hamiltonian operad~\cite{BDK,DotKhor}. 

In the case of linear Poisson brackets, 
all such structures arise from a~pair of compatible Lie brackets.
An algebra $\langle L,[\cdot,\cdot]_1,[\cdot,\cdot]_2,+\rangle$ 
belongs to a~variety $\mathrm{Lie}_2$ of pairs of compatible Lie brackets 
if $\alpha [\cdot,\cdot]_1 + \beta [\cdot,\cdot]_2$ 
is a~Lie bracket for all $\alpha,\beta\in \Bbbk$.

A plenty of examples of compatible Lie brackets is presented in~\cite{GolSokol}, 
the classification results on them see in~\cite{Panasyuk}.
In~\cite{GolSokol2}, it was shown that every pair of compatible Lie brackets endowed
with a common non-degenerate invariant bilinear form produces a~rational solution to the 
classical Yang---Baxter equation.
Free algebras with a~pair of compatible Lie brackets were studied in~\cite{BDK,Leon,Liu}.
Koszulness of the operad corresponding 
to a~variety $\mathrm{Lie}_2$ was proved in~\cite{DotKhor2}.

In~\cite{Khoroshkin}, the operadic (multiplicative) universal enveloping 
associative algebra $U_{\mathrm{Lie}_2}(\mathfrak{g})$ 
of a~given algebra $\mathfrak{g}\in \mathrm{Lie}_2$ 
in the sense of V. Ginzburg and M. Kapranov~\cite{Koszul} was considered, and 
the Poincar\'{e}---Birkhoff---Witt (PBW) property for it was proved.
By the definition, the associative algebra 
$U_{\mathrm{Lie}_2}(\mathfrak{g})$ satisfies the following property:
the category of modules over $\mathfrak{g}$ and 
the category of left modules over $U_{\mathrm{Lie}_2}(\mathfrak{g})$
are equivalent.

We find the Gr\"{o}bner---Shirshov basis of the universal 
enveloping algebra $U_{\mathrm{Lie}_2}(\mathfrak{g}_0)$ 
of an algebra~$\mathfrak{g}_0$, where $\mathfrak{g}_0$ 
denotes the vector space $\mathfrak{g}$ with both zero Lie brackets. 
It allows us, applying the PBW property, to get the linear basis of 
the algebra $U_{\mathrm{Lie}_2}(\mathfrak{g})$. 
As a~corollary, we compute 
the (exponential) growth rate of $U_{\Lie_2}(\mathfrak{g})$ 
when $\mathfrak g$ is finite-dimensional.

\section{Gr\"{o}bner---Shirshov basis for $U_{\mathrm{Lie}_2}(\mathfrak{g}_0)$}

Due to~\cite[Corollary 2.11]{Khoroshkin}, the operadic universal enveloping 
associative algebra of a~given algebra $\mathfrak{g}\in \mathrm{Lie}_2$ 
in the sense of V. Ginzburg and M. Kapranov~\cite{Koszul} equals
\begin{multline} \label{UnivEnvDefRel}
U_{\mathrm{Lie}_2}(\mathfrak{g})
 = \mathrm{As}\langle X\cup X'\mid xy - yx + [x,y]_1,
\,x'y'-y'x'+([x,y]_2)',\\
xy'-y'x+x'y-yx'+([x,y]_1)'+[x,y]_2\rangle,
\end{multline}
where $X$ is a linear basis of~$\mathfrak{g}$, 
$X'$ is a set such that $X\cap X' = \emptyset$ 
and the map $'\colon X\to X'$, $x\to x'$ is a~bijection.
It looks like there are misprints in~\cite[\S6.1]{Khoroshkin} 
while writing the defining relations of $U_{\mathrm{Lie}_2}(\mathfrak{g})$. 

Assuming that $X$ and $X'$ are primitive elements, 
$U_{\mathrm{Lie}_2}(\mathfrak{g})$ has a natural Hopf algebra structure.

In~\cite{Khoroshkin}, the PBW property 
of $U_{\mathrm{Lie}_2}(\mathfrak{g})$ was proved, 
it implies that there exists a~filtration on 
$U_{\mathrm{Lie}_2}(\mathfrak{g})$ such that
$\mathrm{gr}\,U_{\mathrm{Lie}_2}(\mathfrak{g})
 \cong U_{\mathrm{Lie}_2}(\mathfrak{g}_0)$,
where $\mathfrak{g}_0$ is a vector space~$\mathfrak{g}$ 
with trivial products $[\cdot,\cdot]_1$~and~$[\cdot,\cdot]_2$.

Thus, let us study the algebra $U_{\Lie_2}(\mathfrak{g}_0)$.
We may assume that $X = \{f_i\mid i\in I\}$, where $I$ is a well-ordered set.
Denote $X' = \{F_i:= f_i'\mid i\in I\}$.
We define an order on $X\cup X'$ as follows,

\noindent--- $f_i<f_j$ if $i<j$, 

\noindent--- $F_j<F_i$ if $i<j$,

\noindent--- $f_i<F_j$ for all $i,j\in I$.

Let us write down the relations
\begin{gather}
f_j f_i - f_i f_j,\quad i<j, \label{ff} \\
F_i F_j - F_j F_i,\quad i<j, \label{FF} \\
F_i w f_k - f_k F_i w - F_k w f_i + f_i F_k w, \quad w\leq i<k, \label{Ff}
\end{gather}
where $w = 1$ or $w = f_{s_1}\ldots f_{s_t}$ with 
$s_1\leq \ldots \leq s_t$ and by $w\leq i$ we mean that $s_t\leq i$.

We recall the main definitions from the theory of Gr\"{o}bner---Shirshov bases~\cite[\S2.1]{BC}.

Let $(X,<)$ be a well-ordered set and let $X^*$ denote the set of all words in the alphabet~$X$.  
Suppose that $X^*$ is well-ordered, moreover, $u<v$ implies $w_1uw_2<w_1vw_2$ for all $w_1,w_2\in X^*$,
such ordering is called monomial. We will use only deg-lex ordering, in which two words first are compared by the degree and then lexicographically. 
Given a~nonzero element $f$ from the free associative algebra $\As(X)$,
by $\bar{f}$ we mean its leading word.

Given a monomial ordering $<$ on $X^*$ and two monic polynomials $f,g$, we define
two kinds of compositions:

(i) If $w$ is a word such that 
$w = \bar{f} b = a\bar{g}$ for some $a,b\in X^*$ with $|\bar{f}|+|\bar{g}| > |w|$, 
then the polynomial $(f,g)_w := fb - ag$ is called 
the intersection composition of $f$ and $g$ with respect to $w$.

(ii) If $w = \bar{f} = a\bar{g}b$ for some $a, b \in X^*$, then the polynomial $(f,g)_w :=  f - agb$
is called the inclusion composition of $f$ and $g$ with respect to $w$.

Consider $S \subset \As(X)$ such that every $s \in S$ is monic. 
Take $h \in \As(X)$ and $w\in X^*$. 
Then $h$ is called trivial modulo $(S, w)$, denoted by
$h \to 0$ mod $(S, w)$, if 
$h = \sum \alpha_i a_i s_i b_i$, where $\alpha_i \in \Bbbk$, $a_i,b_i\in X^*$, 
and $s_i\in S$ satisfying $a_is_ib_i<w$.

A monic set $S \subset \As(X)$ is called a~Gr\"{o}bner---Shirshov basis in $\As(X)$ 
with respect to the monomial ordering $<$ if every composition of polynomials in $S$ 
is trivial modulo $S$ and the corresponding $w$.

The Composition Diamond lemma for associative algebras implies that
if $S$ is a~Gr\"{o}bner---Shirshov basis in $\As(X)$, then
$\mathrm{Irr}(S) = \{u \in X^*\mid u\neq a\bar{s}b,\, s \in S,\,a,b \in X^*\}$ 
is a linear basis of the algebra $\As\langle X\mid S\rangle$.

{\bf Lemma}.
The set of the relations~\eqref{ff}--\eqref{Ff}
forms a~Gr\"{o}bner---Shirshov basis for $U_{\Lie_2}(\mathfrak{g}_0)$.

{\sc Proof}.
Note that the first series of the defining relations~\eqref{UnivEnvDefRel} of~$U_{\Lie_2}(\mathfrak{g}_0)$ 
coincides with~\eqref{ff}, 
the second one with~\eqref{Ff}, 
and the third series of them multiplied by~$w$ on the right gives exactly~\eqref{Ff}.

The relations~\eqref{ff}--\eqref{Ff} may have 
only compositions of intersection but not inclusion. 

Compositions between~\eqref{ff} and~\eqref{ff} 
as well as between~\eqref{FF} and~\eqref{FF} are trivial.

Let us compute the composition between~\eqref{ff} and~\eqref{Ff}.
Let $w\leq i<k$ and $l<k$, then
\begin{gather*}
F_i w f_k f_l\mathop{\to}\limits_{\eqref{ff}} L:=F_i w f_l f_k, \\
F_i w f_k f_l\mathop{\to}\limits_{\eqref{Ff}}
 R:=f_k F_i w f_l + F_k w f_i f_l - f_i F_k w f_l.
\end{gather*}

If $l\leq i$, then $L-R\mathop{\to}\limits_{\eqref{Ff},\,\eqref{ff}}0$ for $\widetilde{w}= wf_l$.

Let $i<l$, then
\begin{multline*}
L \mathop{\to}\limits_{\eqref{Ff}} 
 f_l F_i w f_k + F_l w f_i f_k - f_i F_l w f_k 
 \mathop{\to}\limits_{\eqref{Ff}} 
 f_l f_k F_i w + f_l F_k w f_i - \underline{f_l f_i F_k w} \\
 + F_l w f_i f_k - f_i f_k F_l w - f_i F_k w f_l + \underline{f_i f_l F_k w};
\end{multline*}
$$
R \mathop{\to}\limits_{\eqref{Ff}} 
 f_k f_l F_i w + f_k F_l w f_i - f_k f_i F_l w + F_k w f_i f_l - f_i F_k w f_l.
$$
Thus,
\begin{multline*}
L-R
 = f_l F_k w f_i + F_l w f_i f_k - f_k F_l w f_i - F_k w f_i f_l \\
 = (f_l F_k w + F_l w f_k - f_k F_l w - F_k w f_l)f_i 
 + F_l w(f_i f_k - f_k f_i) - F_k w (f_i f_l - f_l f_i)
 \mathop{\to}\limits_{\eqref{Ff},\,\eqref{ff}} 0.
\end{multline*}
Here it is important that all involved terms are less than 
the initial word $u = F_i w f_k f_l$.

Now we compute the composition between~\eqref{FF} and~\eqref{Ff}.
Let $w\leq i<k$ and $a<i$. On the one hand, we have
\begin{multline*}
F_a F_i w f_k \mathop{\to}\limits_{\eqref{FF}} 
 F_i F_a w f_k \mathop{\to}\limits_{\eqref{ff}} 
 F_i F_a f_k w \mathop{\to}\limits_{\eqref{Ff}}
 F_i f_k F_a w + F_i F_k f_a w - F_i f_a F_k w \\
 \mathop{\to}\limits_{\eqref{Ff}} 
 L:= f_k F_i F_a w + F_k f_i F_a w - f_i F_k F_a w
 + F_i F_k f_a w - F_i f_a F_k w.
\end{multline*}
On the other hand,
\begin{multline*}
F_a F_i w f_k \mathop{\to}\limits_{\eqref{Ff}}  
 F_a f_k F_i w + F_a F_k w f_i - F_a f_i F_k w \\
 \mathop{\to}\limits_{\eqref{Ff}}
 R:= f_k F_a F_i w + F_k f_a F_i w - \underline{f_a F_k F_i w} 
+ F_a F_k w f_i - f_i F_a F_k w - F_i f_a F_k w + \underline{f_a F_i F_k w}.
\end{multline*}
Thus,
\begin{multline*}
L-R
 = F_k f_i F_a w + F_i F_k f_a w
 - F_k f_a F_i w - F_a F_k w f_i \\
 = F_k(f_i F_a w + F_i f_a w - f_a F_i w - F_a w f_i)
 + (F_i F_k - F_k F_i)f_a w 
 - (F_a F_k - F_k F_a) w f_i
 \mathop{\to}\limits_{\eqref{Ff},\,\eqref{FF}} 0.
\end{multline*}

Finally, note that there are no compositions of intersection between~\eqref{Ff} and~\eqref{Ff}. 
\hfill $\square$

\section{Basis of $U_{\mathrm{Lie}_2}(\mathfrak{g})$}

Denote by $M(X)$ the set of all (ordered) monomials from $\Bbbk[X]$ including~1.
Given $w = f_{j_1}\ldots f_{j_n}\in M(X)\setminus\{1\}$, 
we mean that $j_1\leq \ldots\leq j_n$,
and for $u = F_{k_1}\ldots F_{k_m}\in M(X')\setminus\{1\}$, 
we mean that $k_m\leq \ldots\leq k_1$.
Define 
$$
\lfloor w\rfloor = \max\{j_t\mid t=1,\ldots,n\}, \quad
\lceil u\rceil = \min\{k_t\mid t=1,\ldots,m\},
$$
i.\,e., $\lfloor w\rfloor = j_n$ and $\lceil u\rceil = k_m$. 

Define $L = M(X)\cup L'$, where $L'$ consists of all words
\begin{equation} \label{BaseWord}
w_0 u_1 w_1 u_2 w_2 \ldots u_{s-1}w_{s-1} u_s w_s,
\end{equation}
where 

a) $w_i\in M(X)\setminus\{1\}$, $i=1,\ldots,s-1$, $w_0,w_s\in M(X)$;

b) $u_i\in M(X')\setminus\{1\}$, $i=1,\ldots,s$;

b) $\lfloor w_i\rfloor\leq \lceil u_i\rceil$, $i=1,\ldots,s-1$, 
and $\lfloor w_s\rfloor\leq \lceil u_s\rceil$ or $w_s = 1$.

{\bf Theorem 1}. 
The set $L$ forms a~linear basis of $U_{\Lie_2}(\mathfrak{g})$.

{\sc Proof}.
It follows from Lemma, the Composition Diamond lemma for 
associative algebras~\cite[Theorem 1]{BC}, and the PBW property 
of $U_{\Lie_2}(\mathfrak{g})$~\cite{Khoroshkin}.
\hfill $\square$

Given a pair of compatible Lie brackets~$\mathfrak{g}$,
define $L_n$ as the set of all elements from~$L$ of length~$n$.
Put $r_n = |L_n|$. The growth rate of $\rho(U_{\Lie_2}(\mathfrak{g}))$ is defined as
$\lim\limits_{n\to\infty}\sqrt[n]{r_n}$.
 
Let us show that such limit always exists.
The Fekete's Lemma~\cite[Lemma 1.2.2]{Fekete} 
says that given a sequence $\{a_n\}$, $n\geq1$, 
of real numbers such that $a_{s+t} \leq a_s + a_t$ for all $s,t\in \mathbb{N}$, 
there exists a limit $\lim\limits_{n\to\infty}\frac{a_n}{n}$.
In our case, we have the inequality $r_{s+t} \leq r_s r_t$ for all $s,t\in\mathbb{N}$, 
since every non-empty subword of the basic element from~$L$ lies in~$L$. 
Hence, it remains to apply the Fekete's Lemma for the sequence $\{\ln r_n\}$.

We are able to compute the growth rate of $U_{\Lie_2}(\mathfrak{g})$
when $\mathfrak g$ is finite-dimensional.
Firstly, we do it straightforwardly (Theorem~2).
Secondly, we derive this result from Lemma~6.1~\cite{Khoroshkin} (Remark 1, suggested by the reviwer).
Thirdly, we reprove Theorem~2 with the help of partially commutative algebras and dependence polynomial (Remark~2).

{\bf Theorem 2}.
Let $\mathfrak{g}\in\Lie_2$ and $\dim(\mathfrak g) = m$.
Then the growth rate of $U_{\Lie_2}(\mathfrak{g})$ equals $m+1$.

{\sc Proof}.
Let $X = \{x_1,\ldots,x_m\}$ be a~basis of $\mathfrak{g}$.
Define $O_n$ as a subset of $L_n$ consisting of words starting with $F_i\in X'$.
Put $s_n = |O_n|$, assuming that $s_0 = 1$.

Let us derive the following formula,
\begin{equation} \label{sk:main}
s_k = \sum\limits_{p=1}^{k-1}p\binom{p+m}{p+1}s_{k-1-p} + \binom{k+m-1}{k}.
\end{equation}

A word $v\in O_k$ for $k\geq1$ either consists of only letters from $X'$ 
($\binom{k+m-1}{k}$ choices) or $v$ has the form 
$v = u_1 F_r w_1 v'$, where $r=1,\ldots,m$, $u_1\in M(X')$, 
$w_1\in M(X)\setminus\{1\}$, and $v'\in O_h$ for some $h$ (if $h = 0$,
then $v' = 1$). For the latter case, we initially fix the value of
$p = |u_1|+|w_1|$ and then consider all cases of $t = |w_1|\geq1$.
Hence, we have
\begin{equation} \label{sk:transfer}
s_k = \binom{k+m-1}{k}
 + \sum\limits_{p=1}^{k-1}\left(\sum\limits_{r=1}^m 
\sum\limits_{t=1}^p \binom{t+r-1}{t}\binom{p-t+m-r}{p-t}\right)s_{k-1-p},
\end{equation}
where $\binom{t+r-1}{t}$ is responsible for the choice $w_1\in M(\{x_1,\ldots,x_r\})\setminus\{1\}$
and $\binom{p-t+m-r}{p-t}$ corresponds to the choice of $u_1\in M(\{x_r,\ldots,x_m\}')$.

In~\cite[Theorem 1.4]{Trif} the formula
$$
\sum\limits_{i=0}^m \frac{\binom{m}{i}}{\binom{n+m}{p+i}}
 = \frac{n+m+1}{n+1}\cdot\frac{1}{\binom{n}{p}}
$$
for $0\leq m$ and $0\leq p\leq n$ was stated. 
Applying it, we derive
\begin{multline*}
\sum\limits_{t=0}^p \binom{t+a}{a}\binom{p-t+b}{b} 
 = \sum\limits_{t=0}^p\frac{(t+a)!}{t!a!}\frac{(p-t+b)!}{(p-t)!b!} \\
 = \frac{1}{a!b!}\sum\limits_{t=0}^p\frac{(t+a)!(p-t+b)!}{t!(p-t)!}\frac{p!}{p!}\frac{(p+a+b)!}{(p+a+b)!}
 = \binom{p+a+b}{p,a,b}\sum\limits_{t=0}^p \frac{\binom{p}{t}}{\binom{p+a+b}{t+a}} \\
 = \frac{p+a+b+1}{a+b+1}\binom{p+a+b}{p,a,b}/\binom{a+b}{a}
 = \binom{p+a+b+1}{p}.
\end{multline*}

Substituting this equality for $a = r-1$ and $b = m-r$ in~\eqref{sk:transfer}, we get
\begin{multline*}
\sum\limits_{r=1}^m 
\sum\limits_{t=1}^p \binom{p-t+m-r}{p-t}\binom{t+r-1}{t} \\
 = \sum\limits_{r=1}^m \left(
 \sum\limits_{t=0}^p \binom{p-t+m-r}{p-t}\binom{t+r-1}{t} 
 - \binom{p+m-r}{p}\right) \\
 = \sum\limits_{r=1}^m \left(\binom{p+m}{p} - \binom{p+m-r}{p}\right)
 = m\binom{p+m}{p} - \binom{p+m}{p+1}
 = p\binom{p+m}{p+1},
\end{multline*}
and we have proved the formula~\eqref{sk:main}.

Let us prove that $s_n = m(m+1)^{n-1}$ for any $n\geq1$ by induction on~$n$. The case $n=1$ is trivial, we have $O_1 = \{F_1,\ldots,F_m\}$ and $s_1 = m$.

Suppose that we have proved that $s_n = m(m+1)^{n-1}$ for all $n<k$, where $2\leq k$. By~\eqref{sk:main}, we get
\begin{multline*}
s_k = \sum\limits_{p=1}^{k-2}p\binom{p+m}{p+1}s_{k-1-p} + (k-1)\binom{k+m-1}{k} + \binom{k+m-1}{k} \\
 = \sum\limits_{p=1}^{k-2}p\frac{m(m+1)\ldots(m+p)}{(p+1)!}m(m+1)^{k-2-p}+ \frac{m(m+1)\ldots(m+k-1)}{(k-1)!} \allowdisplaybreaks\\
 = m\left(
 \sum\limits_{p=1}^{k-2}p\frac{((m+1)-1)(m+1)^{k-2-p}}{(p+1)!}\sum\limits_{l=0}^p{ p\brack l}(m+1)^l + \sum\limits_{l=0}^{k-1}\frac{{ k-1\brack l}}{(k-1)!}(m+1)^l
 \right),
\end{multline*}
here ${p\brack l}$ denotes the (unsigned) Stirling number of the first kind. 

The coefficient $A_q$ of $s_k/m$ by~$(m+1)^q$, $q=0,\ldots,k-1$, equals
\begin{multline}
A_q 
 = \sum\limits_{p=k-q-1}^{k-2}\frac{p}{(p+1)!}{ p\brack q+p+1-k}
 - \sum\limits_{p=k-q-2}^{k-2}\frac{p}{(p+1)!}{ p\brack q+p+2-k}
 + \frac{{ k-1\brack q}}{(k-1)!} \\
 = \sum\limits_{p=k-q-1}^{k-2}\left(
 \frac{p}{(p+1)!}{ p\brack q+p+1-k}
 - \frac{(p-1)}{p!}{ p-1 \brack q+p+1-k}
 \right) \\
 +\frac{{ k-1\brack q}}{(k-1)!}-\frac{(k-2){ k-2\brack q}}{(k-1)!}.
\end{multline}
Applying the identity 
${ n+1\brack l} = n{ n\brack l} + { n\brack l-1}$,
we compute for $p>0$
\begin{multline*}
p{ p\brack q+p+1-k} - (p-1)(p+1){ p-1\brack q+p+1-k} \\
 = p{p\brack q+p+1-k}
 -(p+1){p\brack q+p+1-k} 
 +(p+1){ p-1\brack q+p-k} \\
 = -{p\brack q+p+1-k} +(p+1){ p-1\brack q+p-k}.
\end{multline*}
When $p = 0$, we may apply the same formula only adding the summand 
equal to the Kronecker delta $\delta_{q,k-1}$, since ${0\brack 0} = 1$.

Therefore, we have
\begin{multline*}
A_q
 = \frac{{ k-2\brack q-1}}{(k-1)!}
 + \delta_{q,k-1}
 + \sum\limits_{p=k-q-1}^{k-2}\frac{1}{(p+1)!}\left( 
 -{p\brack q+p+1-k} +(p+1){ p-1\brack q+p-k}
 \right) \\
 = \frac{{ k-2\brack q-1}}{(k-1)!} 
 + \delta_{q,k-1} - \frac{{ k-2\brack q-1}}{(k-1)!}
 = \delta_{q,k-1}.
\end{multline*}
Thus, $s_k = m\sum\limits_{q=0}^{k-1}A_q(m+1)^q = m(m+1)^{k-1}$, as required.

Now, we express
\begin{equation} \label{formula:rn}
r_n = \sum\limits_{l=0}^n \binom{l+m-1}{l}s_{n-l},
\end{equation}
this formula counts the leftmost position $(l+1)$ such that 
the $(l+1)$st letter of a given element $z\in L_n$ is from~$X'$. 
The value $l = n$ means that $z\in M(X)$.

By~\eqref{formula:rn}, we get
\begin{equation} \label{formula:rn2}
r_n = \sum\limits_{l=0}^n \binom{l+m-1}{l}s_{n-l}
 = \binom{n+m-1}{n}
 + m\sum\limits_{l=0}^{n-1} \binom{l+m-1}{l}(m+1)^{n-l-1}.
\end{equation}

On the one hand, $r_n \geq s_n = m(m+1)^{n-1}$ and so, $\rho(U_{\Lie_2}(\mathfrak{g}))\geq m+1$. 

On the other hand,
\begin{multline*}
r_n \leq \binom{n+m-1}{n}
 + \sum\limits_{l=0}^{n-1} \binom{l+m-1}{l}(m+1)^{n-l}
 = \sum\limits_{l=0}^n \binom{l+m-1}{l}(m+1)^{n-l} \\
 = \sum\limits_{l=0}^n \frac{m(m+1)\ldots(m+l-1)}{l!}(m+1)^{n-l}
 \leq \sum\limits_{l=0}^n (m+1)^n
 = (n+1)(m+1)^n,
\end{multline*}
which implies $\rho(U_{\Lie_2}(\mathfrak{g}))\leq m+1$.
Hence, $\rho(U_{\Lie_2}(\mathfrak{g})) = m+1$.
\hfill $\square$

{\bf Remark 1} (suggested by the reviwer).
Let $\mathcal{P} = \cup_{n\ge1}\mathcal{P}(n)$ be a symmetric operad.
One may assign to it the generating series
$\chi_{\mathcal{P}}(p_1,p_2,\ldots)
 = \sum_{n\geq1}\chi_{S_n}(\mathcal{P}(n))$,
where $\chi_{S_n}(V) = \sum_{\rho\vdash n}\frac{p_\rho}{z_\rho}\mathrm{tr}_V(\rho)$
is a symmetric function associated with the corresponding $S_n$-character of the symmetric group given in the basis of $p_k = \sum_i x_i^k$.

In Lemma 6.1 of~\cite{Khoroshkin}, it was stated that
\begin{equation} \label{charU0}
\chi_{U_{\Lie_2}^0} = \frac{\sum_{k\geq1} h_k}{1-\sum_{k\geq1} {p_k}}.
\end{equation}
Calculating the generating series for dimensions of $U_{\Lie_2}(\mathfrak{g})$, when $\dim (\mathfrak{g}) = m$, corresponds to the calculation of the plethystic substitution of the polynomial $f(t)=mt$ into~\eqref{charU0}. 
Plethysm with the numerator creates the generating function for the dimensions of the symmetric algebra~$S(\mathfrak{g})$, since $\sum_{k\geq1} h_k$ is the character of the operad Com. 
Plethysm with the denominator gives
\begin{equation} \label{Plethysm}
\frac{1}{1-\sum\limits_{k\geq1} p_k}\circ(mt)
 = \frac{1}{1-\sum\limits_{k\geq1} mt^k}
 = \frac{1}{1-\frac{mt}{1-t}}
 = \frac{1-t}{1-(m+1)t}, 
\end{equation}
so, one gets precisely the generating function for the numbers $s_k = m(m+1)^{k-1}$. 
Altogether, $r_n$ counts all words of the form $w_0 \gamma$,
where $w_0\in M(X)$ and $\gamma\in O_k$ for appropriate~$k$.
Thus, we get the formula~\eqref{formula:rn}, and Theorem~2 follows.   

Moreover, we have confirmed the conjecture of A.~Khoroshkin~\cite[\S6.1]{Khoroshkin} of isomorphism 
$$
U_{\Lie_2^0}(V) \cong S(V)\otimes FL(V)
$$
of Schur functors while projecting on the subspace of words of length~$n$.
Due to Theorem~1, we may present 
$FL(V)$ as the tensor algebra $T(G(V))$, where $G(V)$ is the quotient of
$S(V\oplus V)$ by the ideal $\Span(v'u - u'v\mid u,v\in V)$.
Here we identify $(a,b)\in V\oplus V$ with $a'+b$.

{\bf Remark 2}.
Theorem 2 may be obtained with the help of 
partially commutative algebras.
Given a~graph $G(V,E)$, an associative algebra
$\As(G) = \As\langle V\mid ab = ba,(a,b)\in E\rangle$ 
is called a~partially commutative algebra.
Dependence polynomial~\cite{Fisher1990} of a~graph~$G$ is defined as 
$D(G,x) = 1+\sum\limits_{k=1}^{\omega(G)}(-1)^k c_k(G)x^k$,
where $c_i(G)$ denotes a~number of distinct cliques in~$G$ of the size~$i$
and $\omega(G)$ equals the size of a maximum clique in~$G$.

Consider a graph~$G(V,E)$ with $V = \{f_1,\ldots,f_m\} \cup \{F_1,\ldots,F_m\}$ and 
$$
E = \{(f_i,f_j)\mid i<j\}\cup \{(F_i,F_j)\mid i<j\}
 \cup \{(F_i,f_j)\mid i<j\}.
$$
Note that the set of all basic words from~$L$ coincides with the set of all pairwise distinct words in~$\As(G)$. 
It is known~\cite{Fisher1990} that the generating function of~$\As(G)$
equals to $1/D(G,x)$.
Since 
$$
D(G,x) = \sum\limits_{i=0}^m \binom{m}{i}(-1)^i(i+1)x^i
 = (1-x)^{m-1}(1-(m+1)x),
$$
we conclude that the growth rate of~$U_{\Lie_2}(\mathfrak{g})$ equals $m+1$.
The fact that $\dfrac{1}{(1-x)^{m-1}(1-(m+1)x)}$ is a~generating function of~$U_{\Lie_2}(\mathfrak{g})$ follows from~\eqref{Plethysm}.
Indeed, it is enough to recall that the generating function of the free commutative $m$-generated algebra equals $1/(1-x)^m$.

\section*{Acknowledgements}

The author is grateful to the anonymous referee for pointing out that the defining relations of the universal enveloping of 
a~pair of compatible Lie algebras have to equal~\eqref{UnivEnvDefRel} and for Remark~1.

The author is supported by the grant of the President of the Russian
Federation for young scientists (MK-1241.2021.1.1).

\medskip
\noindent Vsevolod Gubarev \\
Sobolev Institute of Mathematics \\
Acad. Koptyug ave. 4, 630090 Novosibirsk, Russia \\
Novosibirsk State University \\
Pirogova str. 2, 630090 Novosibirsk, Russia \\
e-mail: wsewolod89@gmail.com

\end{document}